\documentclass{article}

\usepackage{amsmath}
\usepackage{url}

\begin{document}

\title{ A New Formula for The Values of Dirichlet Beta Function at Odd Positive
Integers Based on The WZ Method }
\author{ YI JUN CHEN}
\date{}
\maketitle

\begin{abstract}
By using the related results in the WZ theory, a new (as
far as I know) formula for the values of Dirichlet beta function $\beta (s)
= \sum\limits_{n = 1}^{ + \infty } {\frac{( - 1)^{n - 1}}{(2n - 1)^s}}
$ (where $Re(s) > 0$) at odd positive integers was given.
\end{abstract}

\section{Introduction}

It is well known that for Riemann Zeta function $\zeta (s) = \sum\limits_{n
= 1}^{ + \infty } {\frac{1}{n^s}} $ , (where $Re(s) > 1$), Dirichlet Lambda function $\lambda (s) = \sum\limits_{n = 1}^{ + \infty } {\frac{1}{(2n -1)^s}} $, (where $Re(s) > 1$) and Dirichlet Beta function $\beta (s) =
\sum\limits_{n = 1}^{ + \infty } {\frac{( - 1)^{n - 1}}{(2n - 1)^s}}
$, (where $Re(s) > 0$), the following formulas are valid

\begin{equation}
\label{eq1}
\zeta (2n) = \frac{2^{2n - 1}( - 1)^{n - 1}B_{2n} \pi ^{2n}}{(2n)!},
\quad
n \in N
\end{equation}

\begin{equation}
\label{eq2}
\lambda (s) = \frac{2^s - 1}{2^s}\zeta (s),
\end{equation}

\begin{equation}
\label{eq3}
\beta (2n + 1) = \frac{( - 1)^nE_{2n} }{2(2n)!}\left( {\frac{\pi }{2}}
\right)^{2n + 1} \quad ,
\quad
n \in N_0 = N \cup \{0\}
\end{equation}
where $B_n $ is Bernoulli number, $E_n $ is Euler number, which are given by
the following formulas respectively
\[
\frac{x}{e^x - 1} = \sum\limits_{n = 0}^{ + \infty } {\frac{B_n }{n!}x^n} ,
\quad
\frac{2}{e^x + e^{ - x}} = \sum\limits_{n = 0}^{ + \infty } {\frac{E_n
}{n!}x^n} .
\]
And we know that there are no such simple formulas for $\zeta (2n +1)$ and $\beta (2n)$.
In fact, $\zeta (s)$, $\lambda (s)$, $\beta (s)$ are the special cases of
Dirichlet L-function $L_k (s) = \sum\limits_{n = 1}^{ + \infty } {\frac{\chi
_k (n)}{n^s}} $(where $Re(s) > 1$): $\zeta (s) = L_1 (s)$, $\lambda (s) = L_2
(s)$, $\beta (s) = L_{ - 4} (s)$, where $\chi _k $ is Dirichlet character,
see \cite{shimura}. Some general formulas for $L_{ - k} (s)$ and $L_k (s)$ are
given in \cite{shimura}, \cite{zuker1} and \cite{zuker2}. It is worth mentioning that the evaluation of
special values of Dirichlet L-function is an active research field.

In \cite{chen}, based on the framework of the WZ theory (see \cite{wilf}, \cite{Petkovsek}, \cite{zeilberger}), we
obtained a recurrence formula for $\zeta (2l)$

\begin{equation}
\label{eq4}
\varsigma (2l) = \left( {\frac{2^{2l - 1}}{1 - 2^{2l}}} \right)\left\{
{\left[ {\frac{( - 1)^{l + 1}}{4l} + \frac{( - 1)^l}{2}} \right]\frac{\pi
^{2l}}{\Gamma (2l)} + \sum\limits_{j = 1}^{l - 1} {\frac{( - 1)^{l - j}\pi
^{2(l - j)}}{\Gamma (2(l - j) + 1)}\varsigma (2j)} } \right\}
\end{equation}
where $l \in N$, and $\sum\limits_{k = 1}^0 {a(k)} = 0$ is a convention,
which is equivalent to the following classical formula for Bernoulli
polynomial

\begin{equation}
\label{eq5}
B_{2n} \left( {\frac{1}{2}} \right) = \left( {2^{ - 2n + 1} - 1}
\right)B_{2n}
\end{equation}
where $B_n (x)$ is Bernoulli polynomial of order $n$, given by the following
formula:

\[
\frac{te^{tx}}{e^t - 1} = \sum\limits_{n = 0}^{ + \infty } {\frac{B_n
(x)}{n!}t^n} .
\]

In \cite{chen}, it was pointed out that the ideas of getting the recurrence (\ref{eq4}) can be
used to evaluate similar infinite series. In this paper, we obtained a new
(as far as I know) formula for $\beta (2l - 1)$ (where $l \in N$ ) by using the ideas in \cite{chen}. The main
steps are given as follows: we obtained the special values of $\beta (2l -1)$ at $l = 1,2,3$ by using the method in [4] first, then we formulate a
conjecture which give a general formula for $\beta (2l - 1)$, finally we
proved the conjecture. It is worth mentioning that the method of obtaining
$\beta (1)$, $\beta (3)$ and $\beta (5)$ is  consistent with the one of
proving the conjecture. It is also worth mentioning that when I formulate
the conjecture by the special values of $\beta (2l - 1)$, I used
Mathematica 8.0. For the special values of $l$, $l = 1,2,3$ for example , we
can obtain the special values of $\beta (2l - 1)$ by Mathematica 8.0.

The main result in this paper is the following theorem.

\noindent
\textbf{Theorem.} Let $\beta (s) = \sum\limits_{n = 1}^{ + \infty } {\frac{(
- 1)^{n - 1}}{(2n - 1)^s}} $, where $Re(s) > 0$, then we have $\beta (1) =
\frac{\pi }{4}$, $\beta (3) = \frac{\pi ^3}{32}$, $\beta (5) = \frac{5\pi
^5}{1536}$, more generally, for all $l \in N$, we have

\begin{eqnarray}
\beta (2l - 1)& = &\frac{( - 1)^{l + 1}\pi ^{2l - 1}}{2^{2l}}\left[
{\frac{1}{\Gamma (2l - 1)} + 2\sum\limits_{j = 1}^{l - 1} {\frac{( -
1)^j2^{2j}\lambda (2j)}{\Gamma (2l - 2j)\pi ^{2j}}} } \right]\\
 &=& \frac{( - 1)^{l + 1}\pi ^{2l - 1}}{2^{2l}}\left[ {\frac{1}{\Gamma (2l -
1)} + 2\sum\limits_{j = 1}^{l - 1} {\frac{( - 1)^j(2^{2j} - 1)\zeta
(2j)}{\Gamma (2l - 2j)\pi ^{2j}}} } \right]
\end{eqnarray}

\begin{equation}
\label{eq8}
 = \frac{( - 1)^{l + 1}\pi ^{2l - 1}}{2^{2l}}\left[ {\frac{1}{\Gamma (2l -
1)} - \sum\limits_{j = 1}^{l - 1} {\frac{2^{2j}(2^{2j} - 1)B_{2j} }{\Gamma
(2l - 2j)\Gamma (2j + 1)}} } \right].
\end{equation}
By formula (\ref{eq3}) and formula (\ref{eq8}), we obtained the
following formula

\begin{equation}
\label{eq9}
E_{2l} = 1 - \frac{1}{2l + 1}\sum\limits_{j = 1}^l {{{2l + 1} \choose {2j}}} 2^{2j}(2^{2j} - 1)B_{2j}
\end{equation}
where $l \in N_0$, and $\sum\limits_{k = 1}^0 {a(k)} = 0$ is a convention.

Both of formula (6) and formula (7) are similar to formula (4). Because
formula (7) can be obtained by formula (2) and formula (6), and formula (8)
can be obtained by formula (1) and formula (7), we just need to prove formula
(6).

\section{Preparation Lemmas}

To prove the theorem, we need the following lemmas. Lemma 1 can be seen in \cite{chen}, Lemma 2 can be seen in Page 204 of \cite{melzak}, Lemma 3 can be seen in Page 37 of \cite{gradshteyn}, and the proofs of Lemma 4,
Lemma 5 and Lemma 6 will be given below.

\noindent
\textbf{Lemma 1.} For a continuous-discrete WZ pair $(F(x,k),G(x,k))$, that
is, they satisfy the following so-called continuous-discrete WZ
equation

\begin{equation}
\label{eq10}
\frac{\partial F(x,k)}{\partial x} = G(x,k + 1) - G(x,k)
\end{equation}

\noindent
then for all $m,n \in N_0 $, $h,x \in R$, we have

\begin{equation}
\label{eq11}
\sum\limits_{k = m}^n {F(x,k)} - \sum\limits_{k = m}^n {F(h,k)} = \int_h^x
{G(t,n + 1)dt - \int_h^x {G(t,m)} dt} .
\end{equation}

\noindent
\textbf{Lemma 2.} If for all $a,x \in R$, $f(t)$ is integrable on $(a,x)$, then
we have

\begin{equation}
\label{eq12}
\int_a^x {\left( {\int_a^{t_k } { \cdots \left( {\int_a^{t_2 } {f(t_1 )}
dt_1 } \right)} \cdots dt_{k - 1} } \right)} dt_k = \frac{1}{\Gamma
(k)}\int_a^x {(x - t)^{k - 1}f(t)} dt.
\end{equation}

\noindent
\textbf{Lemma 3.} For all $n \in N$, we have

\[
\sum\limits_{k = 1}^n {\cos ((2k - 1)x)} = \frac{1}{2}\sin (2nx)\csc x.
\]

\noindent
\textbf{Lemma 4. }If $k \in N$, we have

\[
\mathop {\lim }\limits_{k \to + \infty } \int_0^{\frac{\pi }{2}} {\frac{\sin
(2kt)}{\sin t}dt = \frac{\pi }{2}} .
\]

\noindent
\textbf{Proof.}
At first, we have

\begin{eqnarray*}
\int_0^{\frac{\pi }{2}} {\frac{\sin (2kt)}{\sin t}dt} & = & \frac{1}{2}\int_0^\pi{\sin \left( {\frac{(2k + 1)t}{2}} \right)\cot \left( {\frac{t}{2}}\right)dt}\\
&& -\: \frac{1}{2}\int_0^\pi {\cos \left( {\frac{(2k + 1)t}{2}}\right)dt}\\
& = & I_1 - I_2 ,
\end{eqnarray*}
where
\[
I_1=\frac{1}{2}\int_0^\pi{\sin \left( {\frac{(2k + 1)t}{2}} \right)\cot \left( {\frac{t}{2}}\right)dt},
\]

\[
I_2=\frac{1}{2}\int_0^\pi {\cos \left( {\frac{(2k + 1)t}{2}}\right)dt}.
\]
Now, let us think about $I_2 $ first, because

\[
\left| {I_2 } \right| = \frac{1}{2}\left| {\int_0^\pi {\cos \left(
{\frac{(2k + 1)t}{2}} \right)dt} } \right| = \frac{1}{2k + 1}\left| {\sin
\left( {\frac{(2k + 1)\pi }{2}} \right)} \right| \le \frac{1}{2k + 1},
\]
we have$\mathop {\lim }\limits_{k \to + \infty } I_2 = 0$.
Next, let us think about $I_1 $. Because $\tan \left( {\frac{t}{4}} \right)$ is continuous on $\left[ {0,\pi }
\right]$, of course, it is integrable on $\left[ {0,\pi } \right]$, by
Riemann-Lebesgue Lemma, we have

\[
\mathop {\lim }\limits_{k \to + \infty } \int_0^\pi {\sin \left( {\frac{(2k+ 1)t}{2}} \right)\tan \left( {\frac{t}{4}} \right)dt}= 0.
\]
By Lemma 4 in \cite{chen} and $\csc \left({\frac{t}{2}}\right)  -\cot \left( {\frac{t}{2}} \right)=\tan \left({\frac{t}{4}} \right),$ we have

\[
\mathop {\lim }\limits_{k \to + \infty } \int_0^\pi {\sin \left( {\frac{(2k
+ 1)t}{2}} \right)\cot \left( {\frac{t}{2}} \right)dt} = \mathop {\lim }\limits_{k \to + \infty }\int_0^\pi {\sin
\left( {\frac{(2k + 1)t}{2}} \right)\csc \left( {\frac{t}{2}} \right)dt } = \pi.
\]
Finally, we have

\begin{eqnarray*}
\mathop {\lim }\limits_{k \to + \infty } \int_0^{\frac{\pi }{2}} {\frac{\sin
\left( {2kt} \right)}{\sin t}dt} &=& \mathop {\lim }\limits_{k \to + \infty }
\left( {I_1 - I_2 } \right) \\
 &=& \frac{1}{2}\mathop {\lim }\limits_{k \to +
\infty } \int_0^\pi {\sin \left( {\frac{(2k + 1)t}{2}} \right)\csc \left(
{\frac{t}{2}} \right)dt} - \mathop {\lim }\limits_{k \to + \infty } I_2 \\
 & =&\frac{\pi }{2} .
\end{eqnarray*}
The proof of Lemma 4 was completed.

\noindent
\textbf{Remarks: 1.} We can also prove

\[
\mathop {\lim }\limits_{k \to + \infty } \int_0^\pi {\sin \left( {\frac{(2k
+ 1)t}{2}} \right)\tan \left( {\frac{t}{4}} \right)dt} = 0
\]

\noindent
by The Second Mean Value Theorem for Integrals. Because for all $t \in \left[
{0,\pi } \right]$, $\left( {\tan \left( {\frac{t}{4}} \right)} \right)^{'} =
\frac{1}{4}\left( {\sec \left( {\frac{t}{4}} \right)} \right)^2 > 0$, we
know that $\tan \left( {\frac{t}{4}} \right)$ is monotone (increasing) on
$\left[ {0,\pi } \right]$, we conclude by The Second Mean Value Theorem for Integrals that there exist $\xi $ on $\left[ {0,\pi } \right]$ , such that

\begin{eqnarray*}
\lefteqn{\left| {\int_0^\pi {\sin \left( {\frac{(2k + 1)t}{2}} \right)\tan \left({\frac{t}{4}} \right)dt} } \right|}\\
 &= &\left| {\tan \left( {\frac{0 + 0}{4}} \right)\int_0^\xi {\sin \left(
{\frac{(2k + 1)t}{t}} \right)} dt + \tan \left( {\frac{\pi - 0}{4}} \right)\int_\xi ^\pi {\sin \left( {\frac{(2k + 1)t}{2}}
\right)dt} } \right|\\
 &\le& 4\tan \left( {\frac{\pi }{4}} \right)\frac{1}{2k + 1}.
\end{eqnarray*}
Finally, we have

\[
\mathop {\lim }\limits_{k \to + \infty } \int_0^\pi {\sin \left( {\frac{(2k
+ 1)t}{2}} \right)\tan \left( {\frac{t}{4}} \right)dt} = 0.
\]

\noindent
\textbf{2.} In fact, this Lemma can be proved in a simpler way by using Lemma 3, the idea here is similar to that in Page 66 of \cite{byerly}. By Lemma 3, for the integrand we have

\[
\frac{\sin(2kt)}{\sin(t)}=2\sum\limits_{k = 1}^n {\cos ((2k - 1)t)}.
\]
Integrate the sum termwise, and note that

 \[
 \int_0^{\frac{\pi }{2}}{\cos((2i-1)t)dt}=\frac{(-1)^{i+1}}{2i-1},
 \]
then we have

\[
\int_0^{\frac{\pi }{2}}{\frac{\sin(2kt)}{\sin(t)}}dt=2\sum\limits_{i=1}^k{\frac{(-1)^{i+1}}{2i-1}},
\]
which tends to $2\arctan(1)=\frac{\pi }{2}$. This completes the proof.

\noindent
\textbf{Lemma 5.} If $k \in N$, for all $s \ge 1$, we have

\[
\mathop {\lim }\limits_{k \to + \infty } \int_0^{\frac{\pi }{2}}
{\frac{t^s\sin (2kt)}{\sin t}dt = 0} .
\]

\noindent
\textbf{Proof.}
Let $f(t) = \left\{ {{\begin{array}{*{20}c}
 {\textstyle{{t^s} \over {\sin t}}} \hfill & {0 < t \le \textstyle{\pi \over
2}} \hfill \\
 0 \hfill & {t = 0} \hfill \\
\end{array} }} \right.$, where $s \ge 1$, then it is easy to prove that
$f(t)$ is monotone (increasing) on $\left[ {0,\frac{\pi }{2}} \right]$ (the
details of proving will be given in the Remarks below) . By The Second Mean Value Theorem for Integrals, we know that there exist $\xi $ on $\left[ {0,\frac{\pi
}{2}} \right]$, such that

\begin{eqnarray*}
\int_0^{\frac{\pi }{2}} {\frac{t^s}{\sin t}\sin (2kt)} dt &=&
\int_0^{\frac{\pi }{2}} {f(t)\sin (2kt)dt}\\
 &=& f(0 + 0)\int_0^\xi {\sin (2kt)dt} + f\left( {\frac{\pi }{2} - 0}
\right)\int_\xi ^{\frac{\pi }{2}} {\sin (2kt)} dt\\
 &=& \left( {\frac{\pi }{2}} \right)^s\int_\xi ^{\frac{\pi }{2}} {\sin (2kt)}dt\\
 &=& \left( {\frac{\pi }{2}} \right)^s\left[ {\frac{1}{2k}\left. {( - \cos
(2kt))} \right|{\begin{array}{*{20}c}
 {\textstyle{\pi \over 2}} \hfill \\
 \xi \hfill \\
\end{array} }} \right].
\end{eqnarray*}
We conclude that

\[
\left| {\int_0^{\frac{\pi }{2}} {\frac{t^s}{\sin t}\sin (2kt)dt} } \right|
\le \left( {\frac{\pi }{2}} \right)^s\frac{1}{k},
\]

\noindent
finally, we have

\[
\mathop {\lim }\limits_{k \to + \infty } \int_0^{\frac{\pi }{2}}
{\frac{t^s}{\sin t}\sin (2kt)} dt = 0.
\]
The proof of Lemma 5 was completed.

\noindent
\textbf{Remarks: 1.} In this remark, we give the proof that $f(t)$ is
monotone (increasing) on $\left[ {0,\frac{\pi }{2}} \right]$. Because $s \ge
1$, it is easy to prove that $f(t)$ is differentiable on $\left[
{0,\frac{\pi }{2}} \right]$, and
$
f^{'}(t) = \left\{ {{\begin{array}{*{20}c}
 {st^{s - 1}\csc t - t^s\cot t\csc t} \hfill & {0 < t \le \textstyle{\pi
\over 2}} \hfill \\
 0 \hfill & {t = 0} \hfill \\
\end{array} }} \right..
$
Let
$
g(t) = s\sin t - t\cos t,
\quad
0 \le t \le \frac{\pi }{2},
$
then we have $g(0) = 0$, $g^{'}(t) = (s - 1)\cos t + t\sin t$. Because $s \ge
1$, for $0 < t \le \frac{\pi }{2}$, we conclude that $(s - 1)\cos t \ge 0$ and $t\sin
t > 0$, that is, for $0 < t \le \frac{\pi }{2}$,  $g^{'}(t) > 0$ . We
conclude that for $0 < t \le \frac{\pi }{2}$,  $g(t) > g(0) = 0$ ,
that is, for $0 < t \le \frac{\pi }{2}$, $f^{'}(t) > f^{'}(0) = 0$. Finally,
because $f(t)$ is continuous on $\left[ {0,\frac{\pi }{2}} \right]$, we
conclude that $f(t)$ is monotone (increasing) on $\left[ {0,\frac{\pi }{2}}
\right]$.

\noindent
\textbf{2.} In fact, the simpler proof of Lemma 5 is the proof  by using Riemann-Lebesgue Lemma directly. Because $f(t) =
\left\{ {{\begin{array}{*{20}c}
 {\textstyle{{t^s} \over {\sin t}}} \hfill & {0 < t \le \textstyle{\pi \over
2}} \hfill \\
 0 \hfill & {t = 0} \hfill \\
\end{array} }} \right.$, where $s \ge 1$, is continuous on $\left[
{0,\frac{\pi }{2}} \right]$, of course, it is integrable on $\left[
{0,\frac{\pi }{2}} \right]$, by  Riemann-Lebesgue Lemma, we
conclude that

\[
\mathop {\lim }\limits_{k \to + \infty } \int_0^{\frac{\pi }{2}} {f(t)\sin
(2kt)} dt = 0.
\]
Finally, because $\int_0^{\frac{\pi }{2}} {\frac{t^s}{\sin t}\sin (2kt)dt =
\int_0^{\frac{\pi }{2}} {f(t)\sin (2kt)dt} } $, we conclude that

\[
\mathop {\lim }\limits_{k \to + \infty } \int_0^{\frac{\pi }{2}}
{\frac{t^s}{\sin t}\sin (2kt)} dt = 0.
\]

\noindent
\textbf{3. }By the proof of this lemma above, we conclude that for $\lambda
\in R$, the following result is still valid for all $s \ge 1$

\[
\mathop {\lim }\limits_{\lambda \to + \infty } \int_0^{\frac{\pi }{2}}
{\frac{t^s}{\sin t}\sin (\lambda t)} dt = 0.
\]

\noindent
\textbf{4.} Similar to Lemma 5 in \cite{chen}, when $s \ge 2$, we can also prove this
lemma by using integration by parts, but when $1 \le s < 2$, this method
can't be used.

\noindent
\textbf{Lemma 6.} Let $I_l (x,k) = \sum\limits_{j = 1}^k {\frac{\sin ((2j -
1)x)}{(2j - 1)^l}} $, $J_l (x,k) = \sum\limits_{j = 1}^k {\frac{\cos ((2j -
1)x)}{(2j - 1)^l}} $, where $k,l \in N$, then we have

\begin{equation}
\label{eq13}
I_{2l + 1} (x,k) = \int_0^x {\left( {\int_0^t { - I_{2l - 1} (t_1 ,k)} dt_1
} \right)dt + \int_0^x {J_{2l} (0,k)dt} } \quad
\end{equation}

\begin{eqnarray}
I_{2l - 1} (x,k)& = & ( - 1)^{l + 1}\int_0^x {\left( {\int_0^{t_{2l - 1} } {
\cdots \left( {\int_0^{t_2 } {\frac{\sin (2kt_1 )}{2\sin t_1 }} dt_1 }
\right)} \cdots dt_{2l - 2} } \right)dt_{2l - 1} }\nonumber\\
&&+\:\sum\limits_{j = 1}^{l- 1} {\frac{( - 1)^{l + j + 1}J_{2j} (0,k)x^{2l - 2j -
1}}{\Gamma (2l - 2j)}}\\
 &=& \frac{( - 1)^{l + 1}}{\Gamma (2l - 1)}\int_0^x {\frac{(x - t)^{2l -2}\sin (2kt)}{2\sin t}dt}\nonumber\\
 &&+\:\sum\limits_{j = 1}^{l - 1} {\frac{( -1)^{l + j + 1}J_{2j} (0,k) x^{2l - 2j - 1}}{\Gamma (2l - 2j)}}
\end{eqnarray}

\noindent
\textbf{Proof.}
Formula (\ref{eq13}) will be proved first. Let

\[
F_{2l + 1} (x,k) = \frac{\sin ((2k - 1)x)}{(2k - 1)^{2l + 1}},
\quad
G_{2l + 1} (x,k) = \sum\limits_{j = 1}^{k - 1} {\frac{\cos ((2j - 1)x)}{(2j
- 1)^{2l}}}
\]

\[
F_{2l} (x,k) = \frac{\cos ((2k - 1)x)}{(2k - 1)^{2l}},
\quad
G_{2l} (x,k) = \sum\limits_{j = 1}^{k - 1} {\frac{ - \sin ((2j - 1)x)}{(2j -
1)^{2l - 1}}}
\]
then it is easy to verify that for $i = 0,1$, $(F_{2l + i} (x,k),G_{2l + i}
(x,k))$ is a WZ pair, by Lemma 1, we have

\[
\sum\limits_{j = 1}^k {F_{2l + i} (x,j)} - \sum\limits_{j = 1}^k {F_{2l + i}
(0,j)} = \int_0^x {G_{2l + i} (t,k + 1)} dt - \int_0^x {G_{2l + i} (t,1)}
dt.
\]
Because (A) for all $j$, $F_{2l} (0,j) = \frac{1}{(2j - 1)^{2l}}$, $F_{2l + 1}
(0,j) = 0$, (B) for all $t$, $G_{2l + i} (t,1) = 0$ by the convention:
$\sum\limits_{k = 1}^0 {a_k } = 0$, (C) it is also easy to verify that

\[
G_{2l + 1} (t,k + 1) = \sum\limits_{j = 1}^k {F_{2l} (x,j)} ,
\quad
G_{2l} (t,k + 1) = - I_{2l - 1} (t,k + 1),
\]
we have

\begin{eqnarray*}
\sum\limits_{j = 1}^k {F_{2l + 1} (x,j)}& =& \int_0^x {G_{2l + 1} (t,k + 1)}dt\\
 &=& \int_0^x {\sum\limits_{j = 1}^k {F_{2l} (t,j)} } dt\\
& =& \int_0^x {\left( {\int_0^t {G_{2l} (t_1 ,k + 1)} dt_1 } \right)} dt +
\int_0^x {J_{2l} (0,k)} dt\\
 &=& \int_0^x {\left( {\int_0^t { - I_{2l - 1} (t_1 ,k + 1)} dt_1 } \right)}
dt + \int_0^x {J_{2l} (0,k)} dt.
\end{eqnarray*}
The proof of formula (\ref{eq13}) is completed.

Now, we will prove formula (14) by mathematical induction. At first, we will prove that when $l = 1,2$, formula (14) is valid. When $l = 1$, let

\[
F_1 (x,k) = \frac{\sin ((2k - 1)x)}{2k - 1},
\quad
G_1 (x,k) = \sum\limits_{j = 1}^{k - 1} {\cos ((2j - 1)x)} ,
\]
it is easy to verify that $(F_1 (x,k),G_1 (x,k))$ is a WZ pair, by
Lemma 1, we have

\[
\sum\limits_{j = 1}^k {F_1 (x,j)} - \sum\limits_{j = 1}^k {F_1 (0,j)} =
\int_0^x {G_1 (t,k + 1)} dt - \int_0^x {G_1 (t,1)} dt
\]
Because (A) for all $j$, $F_1 (0,j) = 0$, (B) for all $t$, $G_1 (t,1) = 0$ by the
convention: $\sum\limits_{k = 1}^0 {a_k } = 0$, by Lemma 3, we have
\begin{eqnarray*}
I_1 (x,k) &=& \sum\limits_{j = 1}^k {F_1 (x,j)}
= \int_0^x {G_1(t,k + 1)} dt\\
 &=&\int_0^x {\sum\limits_{j = 1}^k {\cos ((2j - 1)t)} } dt
 =\frac{1}{2}\int_0^x {\frac{\sin (2kt)}{\sin
t}} dt.
\end{eqnarray*}
By the convention: $\sum\limits_{k = 1}^0 {a_k } = 0$, we have

\[
\sum\limits_{j = 1}^0 {\frac{( - 1)^jJ_{2j} (0,k)}{\Gamma (2 - 2j)}x^{1 -
2j}} = 0.
\]
When $l = 1$, we have

\begin{eqnarray*}
\lefteqn{( - 1)^{l + 1}\int_0^x {\left( {\int_0^{t_{2l - 1} } { \cdots \left(
{\int_0^{t_2 } {\frac{\sin (2kt_1 )}{2\sin t_1 }} dt_1 } \right)} \cdots
dt_{2l - 2} } \right)dt_{2l - 1} }}\\
&=& ( - 1)^{1 + 1}\int_0^x {\frac{\sin(2kt_1 )}{2\sin t_1 }} dt_1
= \frac{1}{2}\int_0^x {\frac{\sin (2kt)}{\sin t}dt}.
\end{eqnarray*}
When $l = 1$, the proof of formula (14) is completed.
When $l = 2$, let

\[
F_3 (x,k) = \frac{\sin ((2k - 1)x)}{(2k - 1)^3},
\quad
G_3 (x,k) = \sum\limits_{j = 1}^{k - 1} {\frac{\cos ((2j - 1)x)}{(2j -
1)^2}}
\]

\[
F_2 (x,k) = \frac{\cos ((2k - 1)x)}{(2k - 1)^2},
\quad
G_2 (x,k) = \sum\limits_{j = 1}^{k - 1} {\frac{ - \sin ((2j - 1)x)}{2j - 1}}
\]

\[
F_1 (x,k) = \frac{ - \sin ((2k - 1)x)}{2k - 1},
\quad
G_1 (x,k) = \sum\limits_{j = 1}^{k - 1} { - \cos ((2j - 1)x)} ,
\]
It is easy to verify that for $i = 1,2,3$, $(F_i (x,k),G_i (x,k))$ is a WZ
pair, by Lemma 1, we have

\[
\sum\limits_{j = 1}^k {F_i (x,j)} - \sum\limits_{j = 1}^k {F_i (0,j)} =
\int_0^x {G_i (t,k + 1)} dt - \int_0^x {G_i (t,1)} dt.
\]
Because (A) for all $j$, $F_3 (0,j) = 0$, $F_2 (0,j) = \frac{1}{(2j - 1)^2}$,
$F_1 (0,j) = 0$, (B) for all $t$, $G_i (t,1) = 0$, $i = 1,2,3$, by the
convention: $\sum\limits_{k = 1}^0 {a_k } = 0$, (C) it is also easy to verify
that

\[
G_{i + 1} (t,k + 1) = \sum\limits_{j = 1}^k {F_i (x,j)} ,
\quad
i = 1,2,
\]
finally, by Lemma 3, we have

\begin{eqnarray*}
I_3 (x,k) &=& \sum\limits_{j = 1}^k {F_3 (x,j)}\\
 &=& \int_0^x {G_3 (t,k + 1)dt}\\
 &=& \int_0^x {\sum\limits_{j = 1}^k {F_2 (t,j)} dt}
\end{eqnarray*}

\begin{eqnarray*}
 &=& \int_0^x {\left( {\int_0^{t_2 } {G_2 (t_1 ,k + 1)} dt_1 } \right)} dt_2 +\int_0^x {J_2 (0,k)} dt\\
 &=& \int_0^x {\left( {\int_0^{t_2 } {\sum\limits_{j = 1}^k {F_1 (t_1 ,j)} }
dt_1 } \right)} dt_2 + \int_0^x {J_2 (0,k)} dt\\
 &=& \int_0^x {\left( {\int_0^{t_3 } {\left( {\int_0^{t_2 } {G_1 (t_1 ,k +1)dt_1 } } \right)dt_2 } } \right)dt_3 + \int_0^x {J_2 (0,k)} dt}\\
 &=& \int_0^x {\left( {\int_0^{t_3 } {\left( {\int_0^{t_2 } {\sum\limits_{j =
1}^k { - \cos ((2j - 1)t_1 )} dt_1 } } \right)dt_2 } } \right)dt_3}\\
&&+\:\int_0^x {J_2 (0,k)} dt\\
 &=& - \int_0^x {\left( {\int_0^{t_3 } {\left( {\int_0^{t_2 } {\frac{\sin
(2kt_1 )}{\sin t_1 }dt_1 } } \right)dt_2 } } \right)dt_3 + J_2 (0,k)x}.
\end{eqnarray*}

\noindent
When $l = 2$, the proof of formula (14) is completed. Next, we will prove that with the assumption that formula (14) is valid for $l$ , then formula (14) is also valid for $l + 1$ . By formula (\ref{eq13}) and the assumption above, we have

\begin{eqnarray*}
\lefteqn{I_{2(l + 1) - 1} (x,k)}\\
&=& I_{2l + 1} (x,k)\\
 &=& \int_0^x {\left( {\int_0^{t_{2l + 1} } { - I_{2l - 1} (t_{2l} ,k)}
dt_{2l} } \right)} dt_{2l + 1} + \int_0^x {J_{2l} (0,k)} dt\\
 &=& ( - 1)^{l + 2}\int_0^x {\left( {\int_0^{t_{2l + 1} } { \cdots \left(
{\int_0^{t_2 } {\frac{\sin (2kt_1 )}{2\sin t_1 }} dt_1 } \right) \cdots }
dt_{2l} } \right)dt_{2l + 1} }\\
 &&+\: \int_0^x {\left( {\int_0^{t_2 } {\sum\limits_{j = 1}^{l - 1} {\frac{( -
1)^{l + j}J_{2j} (0,k)}{\Gamma (2l - 2j)}dt_1 } } } \right)} dt_2 + \int_0^x
{J_{2l} (0,k)} dt\\
 &=& ( - 1)^{l + 2}\int_0^x {\left( {\int_0^{t_{2l + 1} } { \cdots \left(
{\int_0^{t_2 } {\frac{\sin (2kt_1 )}{2\sin t_1 }} dt_1 } \right) \cdots }
dt_{2l} } \right)dt_{2l + 1} }\\
 &&+\: \sum\limits_{j = 1}^{l - 1} {\frac{( - 1)^{l + j}J_{2j} (0,k)}{(2l - 2j +
1)(2l - 2j)\Gamma (2l - 2j)}x^{2l + 2j + 1}} + J_{2l} (0,k)x\\
 &=& ( - 1)^{(l + 1) + 1}\int_0^x {\left( {\int_0^{t_{2l + 1} } { \cdots \left(
{\int_0^{t_2 } {\frac{\sin (2kt_1 )}{2\sin t_1 }} dt_1 } \right) \cdots }
dt_{2l} } \right)dt_{2l + 1} }\\
 &&+\: \sum\limits_{j = 1}^{(l + 1) - 1} {\frac{( - 1)^{(l + 1) + j - 1}J_{2j}
(0,k)}{\Gamma (2(l + 1) - 2j)}x^{2(l + 1) - 2j - 1}} .
\end{eqnarray*}

\noindent
We proved that when formula (14) is valid for $l$ , formula (14) is also valid for $l + 1$ . Finally, by the principle of mathematical induction, we proved that formula (14) is valid for all $l \in n$. By
formula (\ref{eq13}) and Lemma 2, it is easy to prove formula (15). The proof
of Lemma 6 is completed.

\section{Proof of Theorem}

As mentioned above, we just need to prove formula (6). By Lemma 6,
for all $l \in N$, we have

\[
\beta (2l - 1) = \mathop {\lim }\limits_{k \to + \infty } I_{2l - 1} \left(
{\frac{\pi }{2},k} \right),
\quad
\lambda (2l) = \mathop {\lim }\limits_{k \to + \infty } J_{2l} (0,k).
\]
By formula (15), we have

\begin{eqnarray*}
I_{2l - 1} (x,k) &=& \frac{( - 1)^{l + 1}}{\Gamma (2l - 1)}\int_0^x {(x -
t)^{2l - 2}\frac{\sin (2kt)}{2\sin t}dt}\\
&&+\: \sum\limits_{j = 1}^{l - 1}
{\frac{( - 1)^{l + j + 1}J_{2j} (0,k)}{\Gamma (2l - 2j)}x^{2l - 2j - 1}}\\
 &=& I_{2l - 1,1} (x,k) + I_{2l - 1,2} (x,k),
\end{eqnarray*}
where

\[
I_{2l - 1,1} (x,k)=\frac{( - 1)^{l + 1}}{\Gamma (2l - 1)}\int_0^x {(x -
t)^{2l - 2}\frac{\sin (2kt)}{2\sin t}dt},
\]

\[
I_{2l - 1,2} (x,k)=\sum\limits_{j = 1}^{l - 1}
{\frac{( - 1)^{l + j + 1}J_{2j} (0,k)}{\Gamma (2l - 2j)}x^{2l - 2j - 1}}.
\]
It is easy to see that

\begin{eqnarray*}
\mathop {\lim }\limits_{k \to + \infty } I_{2l - 1,2} \left( {\frac{\pi
}{2},k} \right) &=& \sum\limits_{j = 1}^{l - 1} {\frac{( - 1)^{l + j +
1}}{\Gamma (2l - 2j)}\left( {\frac{\pi }{2}} \right)^{2l - 2j - 1}\mathop
{\lim }\limits_{k \to + \infty } J_{2l} (0,k)}\\
 &=& \sum\limits_{j = 1}^{l - 1} {\frac{( - 1)^{l + j + 1}}{\Gamma (2l -
2j)}\left( {\frac{\pi }{2}} \right)^{2l - 2j - 1}\lambda (2j)} .
\end{eqnarray*}
By Lemma 4 and Lemma 5, we have

\begin{eqnarray*}
\lefteqn{\mathop {\lim }\limits_{k \to + \infty } I_{2l - 1,1} \left( {\frac{\pi
}{2},k} \right)}\\
 &=& \mathop {\lim }\limits_{k \to + \infty } \frac{( - 1)^{l +
1}}{\Gamma (2l - 1)}\int_0^{\frac{\pi }{2}} {\left( {\frac{\pi }{2} - t}
\right)^{2l - 2}\frac{\sin (2kt)}{2\sin t}dt}\\
 &=& \frac{( - 1)^{l + 1}}{2\Gamma (2l - 1)}\sum\limits_{j = 1}^{2l - 2}
{{{2l - 2} \choose {j}}\left( {\frac{\pi }{2}} \right)^{2l - 2 - j}} \mathop
{\lim }\limits_{k \to + \infty } \int_0^{\frac{\pi }{2}} {\frac{t^j\sin
(2kt)}{\sin t}} dt
\end{eqnarray*}

\begin{eqnarray*}
&&+\: \frac{( - 1)^{l + 1}}{2\Gamma (2l - 1)}\left( {\frac{\pi }{2}}
\right)^{2l - 2}\mathop {\lim }\limits_{k \to + \infty } \int_0^{\frac{\pi
}{2}} {\frac{\sin (2kt)}{\sin t}} dt\\
&=& \frac{( - 1)^{l + 1}\pi ^{2l - 1}}{2^{2l}\Gamma (2l - 1)}.
\end{eqnarray*}
Finally, we have

\begin{eqnarray*}
\lefteqn{\beta (2l - 1)}\\
& =& \mathop {\lim }\limits_{k \to + \infty } I_{2l - 1} \left(
{\frac{\pi }{2},k} \right)\\
 &=& \mathop {\lim }\limits_{k \to + \infty } I_{2l - 1,1} \left( {\frac{\pi
}{2},k} \right) + \mathop {\lim }\limits_{k \to + \infty } I_{2l - 1,2}
\left( {\frac{\pi }{2},k} \right)\\
 &=& \frac{( - 1)^{l + 1}\pi ^{2l - 1}}{2^{2l}}\left[ {\frac{1}{\Gamma (2l -
1)} + 2\sum\limits_{j = 1}^{l - 1} {\frac{( - 1)^j2^{2j}\lambda (2j)}{\Gamma
(2l - 2j)\pi ^{2j}}} } \right].
\end{eqnarray*}
The proof of Theorem is completed.


\bigskip

\noindent\textit{School of Mathematical Science,
South China Normal University, Guangzhou
510631,P.R.China\\
chenyijun73@yahoo.com.cn}

\end{document}